\documentclass{amsart}


\usepackage{amsmath,amssymb, amstext,amsthm,amsfonts}
\usepackage[dvips]{graphicx}
\usepackage{epic,eepic}

\title[Lyapunov stability and sectional-hyperbolicity for flows]
      {Lyapunov stability and Sectional-hyperbolicity for higher-dimensional flows}

\author[A. Arbieto, C. A. Morales, B. Santiago]{A. Arbieto, C. A. Morales, B. Santiago}
\address{Instituto de Matem\'atica, Universidade Federal do Rio de Janeiro, P. O. Box 68530, 21945-970 Rio
de Janeiro, Brazil.}
\email{arbieto@im.ufrj.br, morales@impa.br, brsantiago777@gmail.com}

\thanks{Partially supported by CNPq, FAPERJ and PRONEX/DS from Brazil.}

\newtheorem{theorem}{Theorem}

\newtheorem*{main}{Theorem}
\newtheorem*{mainc}{Corollary}
\newtheorem{lemma}[theorem]{Lemma}

\newtheorem{definition}[theorem]{Definition}

\newcommand{\cl}{\operatorname{Cl}}

\newcommand{\per}{\operatorname{Per}}

\begin{document}

\begin{abstract}
We study $C^1$-generic vector fields on closed manifolds
without points accumulated by periodic orbits of different indices and prove
that they
exhibit finitely many sinks and sectional-hyperbolic transitive Lyapunov stable sets
with residual basin of attraction.
This represents a partial positive answer to conjectures in \cite{am},
the Palis conjecture \cite{pa} and extend the Araujo's thesis to higher dimensions \cite{a}.
\end{abstract}

\maketitle

\section{Introduction}

\noindent
Dynamical systems (i.e. vector fields or diffeomorphisms)
on closed manifolds and, specifically, the $C^1$ generic ones, have been studied
during these last fifty years or so.
In fact, Pugh proved in the early sixties \cite{p}
that such systems display dense closed orbits in their nonwandering set. Moreover,
Ma\~n\'e \cite{M} proved that a $C^1$ generic surface diffeomorphism is Axiom A or exhibits
infinitely many attracting periodic orbits up to time reversing whereas
Araujo in his thesis \cite{a} proved that
these diffeomorphisms have either infinitely many attracting periodic orbits
or finitely many hyperbolic attractors with full Lebesgue measure basins of attraction.
On the other hand, Hayashi and Ma\~n\'e proved the celebrated Palis-Smale's {\em $C^1$ stability conjecture} \cite{pm}
that all $C^1$ structural stable systems are Axiom A \cite{h}, \cite{mane} and
Ma\~n\'e \cite{M} initiated the study of what today's we call
{\em star systems}, i.e., dynamical systems which cannot be $C^1$-approximated by ones exhibiting
nonhyperbolic closed orbits. By noting that star flows may not be Axiom A
(e.g. the {\em geometric Lorenz attractor} \cite{abs}, \cite{gu}, \cite{GW}),
he asked if, on the contrary, all star diffeomorphisms on closed manifolds are Axiom A.
Such a problem was solved in positive by Aoki and Hayashi \cite{ao}, \cite{h0} and
inspired Gan and Wen \cite{gw} to identify the presence of singularities
in the preperiodic set (c.f. \cite{w}) as the sole obstruction for
a star flow to be Axiom A. In particular, \cite{gw} proved that all nonsingular star flows on closed manifolds
are Axiom A. This solved in positive a conjecture by Liao and Ma\~n\'e.

Meanwhile
\cite{mpp} introduced the notion of
{\em singular-Axiom A flow} inspiried on both the Axiom A flows and the geometric Lorenz attractor.
Based on techniques introduced by Hayashi and Ma\~n\'e for the solution of the stability conjecture,
it was proved in \cite{mpa1} that a $C^1$ generic vector field on a closed $3$-manifold either is singular-Axiom A
or exhibits infinitely many attracting periodic orbits up to flow reversing.
This result motivated the question whether analogous result holds for $C^1$ generic vector fields
in higher dimensional manifolds but negative results were then obtained.
These results motivated the notion of {\em sectional-Axiom A flow} \cite{memo}
as a natural substitute of the singular-Axiom A flows in higher-dimensions.
Unfortunately, results like \cite{mpa1} with the term sectional-Axiom A in place of singular-Axiom A
are not longer true. Instead, the first author conjectured in \cite{am} that a $C^1$ generic
star flow on a closed manifold is sectional-Axiom A. If this conjecture were true, then it would be also true that
all $C^1$-generic vector fields
without points accumulated by hyperbolic periodic orbits of different Morse indices are
sectional-Axiom A.
It was this last assertion what was proved in \cite{am} but when the singularities accumulated by periodic orbits
have Morse index $1$ or $n-1$.

In this paper we prove that all $C^1$ generic vector fields
without points accumulated by periodic orbits of different indices on a closed manifold
are essentially sectional-Axiom A. By this we mean that they
come equipped with
finitely many sinks and sectional-hyperbolic transitive Lyapunov stable sets with residual basin of attraction.
This result (which applies to the star flows with spectral decomposition) represents a partial positive answer
to the aforementioned conjectures \cite{am},
the Palis conjecture \cite{pa} and extend the Araujo's thesis to higher dimensions.
It is also related to \cite{cp} where it was proved that every
$C^1$ diffeomorphism of a closed manifold is approximated by another diffeomorphism with a homoclinic
or heteroclinic bifurcation or by one which is
essentially Axiom A (i.e. with finitely many hyperbolic attractors with open and dense basin of attraction).
Let us state our result in a precise way.

In what follows $M$ will denote a {\em closed $n$-manifold}, i.e.,
a compact connected boundaryless Riemannian manifold of dimension $n\geq 3$.
The space of $C^1$ vector fields in $M$ will be denoted by $\mathcal{X}^1$.
If $X\in\mathcal{X}^1$ we
denote by $X_t$ the flow generated by $X$ in $M$. A {\em periodic orbit} (resp. {\em singularity}) of $X$
is the orbit $\{X_t(p):t\in\mathbb{R}\}$ of a point $p\in M$
satisfying $X_T(p)=p$ for some minimal $T>0$ (resp. a zero of $X$).
By a {\em closed orbit} we mean a periodic orbit or a singularity.
Denote by $Sing_X(\Lambda)$ the set of singularities of $X$ in a subset $\Lambda \subset M$.

Given $p\in M$ we define the {\em $\omega$-limit set}
$$
\omega(p)=\left\{x\in M:
x=\lim_{n\to\infty}X_{t_n}(p)
\mbox{ for some sequence }t_n\to\infty\right\}.
$$
A subset $\Lambda\subset M$ is
{\em invariant} if $X_t(\Lambda)=\Lambda$ for all $t\in \mathbb{R}$;
{\em nontrivial} if it does not reduces to a single closed orbit;
a {\em limit cycle} if $\Lambda=\omega(x)$ for some $x\in M$;
{\em transitive} if $\Lambda=\omega(p)$ for some $p\in \Lambda$ and 
{\em Lyapunov stable}
if for every neighborhood $U$ of it there is a neighborhood $\Lambda\subset W\subset U$ such that
$X_t(W)\subset U$ for all $t\geq 0$.
Moreover, we say that $\Lambda$ has {\em dense closed} (resp. {\em periodic) orbits}
if the closed (resp. periodic) orbits of $X$ in $\Lambda$ are dense in $\Lambda$.
We also define the {\em basin of attraction} of $\Lambda$ by
$$
W^s(\Lambda)=\{x\in M:\omega(x)\subset \Lambda\}.
$$

A transitive set $\Lambda$ will be called {\em attractor}
if it exhibits a neighborhood $U$ such that
$
\Lambda=\bigcap_{t\geq 0}X_t(U).
$
On the other hand, a compact invariant set $\Lambda$ is {\em $C^1$ robustly transitive} if
there is a compact neighborhood $U$ of $\Lambda$ with
$\Lambda=\bigcap_{t\in \mathbb{R}}X_t(U)$ such that
$\Lambda(Y)=\bigcap_{t\in\mathbb{R}}Y_t(U)$
is a nontrivial transitive set of $Y$ for every vector field $Y$ that is $C^1$ close to $X$
($\Lambda(Y)$ is often referred to as the natural continuation of $\Lambda$).

Denote by $\|\cdot\|$ and $m(\cdot)$ the norm and the minimal norm
induced by the Riemannian metric and by
$Det(\cdot)$ the jacobian operation.
We say that $\Lambda$ is {\em hyperbolic}
if there are a continuous invariant tangent bundle decomposition
$$
T_\Lambda M=\hat{E}^s_\Lambda\oplus \hat{E}^X_\Lambda\oplus \hat{E}^u_\Lambda
$$
and positive constants $K,\lambda$
such that $\hat{E}^X_\Lambda$ is the subbundle generated by $X$,
$$
\|DX_t(x)/\hat{E}^s_x\|\leq Ke^{-\lambda t}
\quad \mbox{ and }\quad m(DX_t(x)/\hat{E}^u_x)\geq K^{-1}e^{\lambda t},
$$
for all $x\in \Lambda$ and $t\geq 0$.
A closed orbit is hyperbolic if it does as a compact invariant set.
We define the {\em Morse index} $I(O)$ of a hyperbolic closed orbit $O$ by $I(O)=dim(E^s_x)$
for some (and hence for all) $x\in O$.
In case $O$ is a singularity $\sigma$ we write
$I(\sigma)$ instead of $I(\{\sigma\})$.
A {\em sink} will be a hyperbolic closed orbit of maximal Morse index.

Given an invariant splitting
$T_\Lambda M=E_\Lambda\oplus F_\Lambda$
over an invariant set $\Lambda$ of a vector field $X$
we say that the subbundle $E_\Lambda$ {\em dominates} $F_\Lambda$ if there are positive constants
$K,\lambda$ such that
$$
\frac{\|DX_t(x)/E_x\|}{m(DX_t(x)/F_x)}\leq Ke^{-\lambda t},
\quad\quad\forall x\in \Lambda \mbox{ and }t\geq 0.
$$
(In such a case we say that $T_\Lambda M=E_\Lambda\oplus F_\Lambda$ is a {\em dominated splitting}).

We say that $\Lambda$ is {\em partially hyperbolic} if
it has a dominated splitting $T_\Lambda M=E^s_\Lambda\oplus E^c_\Lambda$
whose dominating subbundle $E^s_\Lambda$ is {\em contracting},
namely,
$$
\|DX_t(x)/E^s_x\|\leq Ke^{-\lambda t},
\quad\quad\forall x\in \Lambda \mbox{ and }t\geq 0.
$$
Moreover, we call the central subbundle $E^c_\Lambda$
{\em sectionally expanding} if
$$
dim(E^c_x)\geq 2 \quad\mbox{ and }\quad
|Det(DX_t(x)/L_x)|\geq K^{-1}e^{\lambda t},
\quad\quad\forall x\in \Lambda \mbox{ and }t\geq 0
$$
and all two-dimensional subspace
$L_x$ of $E^c_x$.

We call {\em sectional-hyperbolic}
any partially hyperbolic set whose singularities (if any) are hyperbolic
and whose central subbundle is sectionally expanding \cite{memo}.

Notice that $\mathcal{X}^1$ is a Baire space if equipped with the standard $C^1$ topology.
We shall use consistently the expression {\em residual subset}
which indicates a certain subset in a metric
space which is a countable intersection of open and dense subsets.
A fundamental property of the set of residual subsets is that
it is closed under countable intersection.
This property will be used implicitely along the proof of our theorem.
We also use the customary expression
{\em $C^1$-generic vector field} meaning {\em for every vector field in a residual subset of $\mathcal{X}^1$}.

With these definitions we can state our main result.

\begin{main}
\label{thA}
Let $X\in\mathcal{X}^1$ be a $C^1$-generic vector field
without points accumulated by hyperbolic periodic orbits of different Morse indices.
Then, $X$ has finitely many sinks and sectional-hyperbolic transitive Lyapunov stable sets
with residual basin of attraction.
\end{main}

The proof will use
some recent results like \cite{glw}, \cite{gwz}, \cite{m}, \cite{mpa2}.
It would be nice to obtain attractors instead of transitive Lyapunov stable sets in this theorem.
Unfortunately, as asked in \cite{cm},
it is unkown whether a sectional-hyperbolic transitive Lyapunov stable set is an attractor (even generically).
Let us present a short application of our result.

We say that $X\in\mathcal{X}^1$ is a {\em star flow}
if there is a neighborhood $\mathcal{U}$ of $X$ such that every
closed orbit of every $Y\in \mathcal{U}$ is hyperbolic.
Recall that the {\em nonwandering set}
of $X$ is the set of points $p\in M$ such that
for every neighborhood $U$ of $p$ and every $T>0$ there is
$t>T$ such that $X_t(U)\cap U\neq\emptyset$.
We say that $X$ has {\em spectral decomposition} if $\Omega(X)$ splits into finitely many disjoint transitive sets.
Moreover, we say that $X$ is
a {\em sectional-Axiom A flow} if there is a finite disjoint union
$\Omega(X)=\Omega_1\cup \cdots \cup\Omega_k$
formed by transitive sets with dense closed orbits
$\Omega_1,\cdots, \Omega_k$ such that, for all $1\leq i\leq k$,
$\Omega_i$ is either a hyperbolic set for $X$ or
a sectional-hyperbolic set for $X$ or a sectional-hyperbolic set for $-X$.
Clearly a sectional-Axiom A flow has a spectral decomposition but the converse is not necessarily true.

As already mentioned, the first author conjectured in \cite{am} that all $C^1$-generic star flows on closed manifolds
are sectional-Axiom A. A support for this conjecture is given below.
Its proof follows from the Theorem and Lemma 16 in \cite{am}.

\begin{mainc}
\label{the-coro}
A $C^1$-generic star flow with spectral decomposition
has finitely many sectional-hyperbolic transitive Lyapunov stable sets with residual basin of attraction.
\end{mainc}

\section{Proof}

\noindent
Previously we state some basic results.
The first one is the main result in \cite{mpa2}.

\begin{lemma}
 \label{l1}
For every $C^1$-generic vector field $X\in \mathcal{X}^1$ there is a residual subset $R_X$ of $M$
such that $\omega(x)$ is a Lyapunov stable set, $\forall x\in R_X$.
\end{lemma}

With the same methods as in \cite{cmp} and \cite{mpa2} it is possible to prove the following variation of this lemma.
We shall use the standard stable and unstable manifold operations $W^s(\cdot),W^u(\cdot)$
(c.f. \cite{hps}).

\begin{lemma}
\label{l2}
For every $C^1$-generic vector field $X\in \mathcal{X}^1$ 
and every hyperbolic closed orbit $O$ of $X$
the set $\{x\in W^u(O)\setminus O:\omega(x)$ is Lyapunov stable$\}$ is
nonempty (it is indeed residual in $W^u(O)$).
\end{lemma}

A straighforward extension to higher dimensions of the three-dimensional arguments
in \cite{m} allows us to prove the following lemma.

\begin{lemma}
 \label{l5}
A sectional-hyperbolic set $\Lambda$ of $X\in \mathcal{X}^1$ contains finitely many attractors, i.e., the collection
$
\{A\subset \Lambda: A\mbox{ is an attractor of }X\}
$
is finite.
\end{lemma}

The following concept comes from \cite{glw}.

\begin{definition}
We say that a compact invariant set $\Lambda$  of $X\in \mathcal{X}^1$ 
has a definite index $0\leq Ind(\Lambda)\leq n-1$ if there are
a neighborhood $\mathcal{U}$ of $X$ in $\mathcal{X}^1$ and
a neighborhood $U$ of $\Lambda$ in $M$ such that
$I(O)=Ind(\Lambda)$ for
every hyperbolic periodic orbit $O\subset U$ of every vector field $Y\in \mathcal{U}$.
In such a case we say that $\Lambda$ is {\em strongly homogeneous (of index $Ind(\Lambda)$)}.
\end{definition}

The importance of the strongly homogeneous property is given by the following result proved in \cite{glw}:
If a strongly homogeneous sets $\Lambda$ with singularities
(all hyperbolic) of $X\in \mathcal{X}^1$ is $C^1$ robustly transitive, then it is
partially hyperbolic for either $X$ or $-X$ depending on whether
\begin{equation}
\label{eql1}
I(\sigma)>Ind(\Lambda),
\quad\quad\forall \sigma\in Sing_X(\Lambda)
\end{equation}
or
\begin{equation}
 \label{eql11}
I(\sigma)\leq Ind(\Lambda),
\quad\quad\forall \sigma\in Sing_X(\Lambda)
\end{equation}
holds.
This result was completed in \cite{memo} by proving that all such sets are in fact
sectional-hyperbolic for either $X$ or $-X$ depending on whether (\ref{eql1}) or (\ref{eql11}) holds.
Another proof of this completion can be found in \cite{gwz}.

On the other hand, \cite{am} observed that the completion in \cite{memo} (or \cite{gwz}) is also valid for
transitive sets with singularities (all hyperbolic of Morse index $1$ or $n-1$)
as soon as $n\geq 4$ and $1\leq Ind(\Lambda)\leq n-2$.
The proof is the same as \cite{glw} and \cite{memo}
but with the preperiodic set playing the role of the
natural continuation of a $C^1$ robustly transitive set.

Now we observe that such a completion is still valid for
limit cycles or when the periodic orbits are dense.
In other words, we have the following result.

\begin{lemma}
\label{l3}
If a strongly homogeneous set $\Lambda$
with singularities (all hyperbolic) of $X\in \mathcal{X}^1$ satisfying
$1\leq Ind(\Lambda)\leq n-2$
is a limit cycle or has dense periodic orbits,
then it is sectional-hyperbolic for either $X$ or $-X$ depending on whether
(\ref{eql1}) or (\ref{eql11}) holds.
\end{lemma}

This lemma motivates the problem whether a strongly homogeneous set with hyperbolic singularities
which is a limit cycle or has dense periodic orbits satisfies either (\ref{eql1}) or (\ref{eql11}).
For instance, Lemma 3.3 of \cite{gwz} proved this is the case for all $C^1$ robustly transitive strongly homogeneous sets.
Similarly for
strongly homogeneous limit cycles with singularities (all hyperbolic of Morse index $1$ or $n-1$)
satisfying $n\geq 4$ and $1\leq Ind(\Lambda)\leq n-2$ (e.g. Proposition 7 in \cite{am}).
Consequently, all such sets
are sectional-hyperbolic for either $X$ or $-X$.
See Theorem A in \cite{gwz} and Corollary 8 in \cite{am} respectively.

Unfortunately,
(\ref{eql1}) or (\ref{eql11}) need not be valid for general
strongly homogeneous sets with dense periodic orbits
even if $1\leq Ind(\Lambda)\leq n-2$.
A counterexample is the nonwandering set of the vector field in $S^3$ obtained by gluing
a Lorenz attractor and a Lorenz repeller as in p. 1576 of \cite{mpa1}.
Despite, it is still possible to analyze the singularities of
a strongly homogeneous set with dense periodic orbits
even if (\ref{eql1}) or (\ref{eql11}) does not hold.
For instance, adapting the proof of Lemma 2.2 in \cite{gwz} 
(or the sequence of lemmas 4.1, 4.2 and 4.3 in \cite{glw}) we obtain the following result.

\begin{lemma}
 \label{l4}
If $\Lambda$ is a strongly homogeneous set with singularities
(all hyperbolic) and dense periodic orbits of $X\in \mathcal{X}^1$, then every
$\sigma\in Sing_X(\Lambda)$ satisfying $I(\sigma)\leq Ind(\Lambda)$
exhibits a dominated splitting $\hat{E}^u_\sigma=E^{uu}_\sigma\oplus E^c_\sigma$
with $dim(E^{uu}_\sigma)=n-Ind(\Lambda)-1$ over $\sigma$ such that the strong unstable manifold $W^{uu}(\sigma)$
tangent to $E^{uu}_\sigma$ at $\sigma$ (c.f. \cite{hps}) satisfies
$\Lambda\cap W^{uu}(\sigma)=\{\sigma\}$.
\end{lemma}

Now we can prove our result.

\begin{proof}[Proof of the Theorem]
Let $X\in \mathcal{X}^1$ be a $C^1$-generic vector field without points accumulated by
hyperbolic periodic orbits of different Morse indices.
By \cite{am}, since $X$ is $C^1$ generic, it follows that
if $\per_i(X)$ denotes the union of the periodic orbits with Morse index $i$,
then the closure $\cl(\per_i(X))$ is
strongly homogeneous of index $Ind(\cl(\per_i(X)))=i$,
$\forall 0\leq i\leq n-1$.
Moreover, $X$ is a star flow and so
it has finitely many singularities and also finitely many
periodic orbits of Morse index $0$ or $n-1$ (c.f. \cite{li}, \cite{pl}).

Let us prove that $\omega(x)$ is sectional-hyperbolic for all $x\in R_X$
where $R_X\subset M$ is the residual subset in Lemma \ref{l1}.
We can assume that $\omega(x)$ is nontrivial and has singularities for, otherwise,
$\omega(x)$ is hyperbolic by Theorem B in \cite{gw} and the Pugh's closing lemma \cite{p}.

Since $X$ is $C^1$ generic we can further assume that $\omega(x)\subset \cl(\per_i(X))$ for some
$0\leq i\leq n-1$ by the closing lemma once more.
Since $X$ has finitely many singularities and periodic orbits of Morse index
$0$ or $n-1$ we have $1\leq i\leq n-2$
(otherwise $\omega(x)$ will be reduced to a singleton which is absurd).
Since $\cl(\per_i(X))$ is strongly homogeneous of index $i$ we have that $\omega(x)$ also does
so $1\leq Ind(\omega(x))\leq n-2$.
Then, since $\omega(x)$ is a limit cycle,
we only need to prove by Lemma \ref{l3} that (\ref{eql1}) holds for $\Lambda=\omega(x)$.
To prove it we proceed as in Corollary B in \cite{glw}, namely,
suppose by contradiction that (\ref{eql1}) does not hold. Then,
there is $\sigma\in Sing_X(\Lambda)$ such that
$I(\sigma)\leq Ind(\omega(x))$.
Since $\omega(x)\subset \cl(\per_i(X))$ and $\cl(\per_i(X))$
is a strongly homogeneous set with singularities,
all hyperbolic, in $\Omega(X)$
we have by Lemma \ref{l4} that
there is a dominated splitting $\hat{E}^u_\sigma=E^{uu}_\sigma\oplus E^c_\sigma$
for which the associated strong unstable manifold $W^{uu}(\sigma)$
satisfies $\cl(\per_i(X))\cap W^{uu}(\sigma)=\{\sigma\}$.
However $W^{uu}(\sigma)\subset \omega(x)$ since $\sigma\in\omega(x)$ and
$\omega(x)$ is Lyapunov stable. As $\omega(x)\subset\cl(\per_i(X))$ we conclude that
$W^{uu}(\sigma)=\{\sigma\}$ so
$\dim(E^{uu}_\sigma)=0$.
But $\dim(E^{uu}_\sigma)=n-i-1$ by Lemma \ref{l4}
so $\dim(E^{uu}_\sigma)\geq n-n+2-1=1$ a contradiction.
We conclude that (\ref{eql1}) holds so $\omega(x)$ is sectional-hyperbolic for all $x\in R_X$.

Next we prove that $\omega(x)$ is transitive for $x\in R_X$.
If $\omega(x)$ has no singularities, then it is hyperbolic
and so a hyperbolic attractor of $X$.
Otherwise, there is $\sigma\in Sing_X(\Lambda)$.
By Lemma \ref{l2} we can select  $y\in W^u(\sigma)\setminus\{\sigma\}$ with Lyapunov stable $\omega$-limit set.
On the other hand, $\omega(x)$ is Lyapunov stable and so
$W^u(\sigma)\subset \omega(x)$. Then,
we obtain $y\in \omega(x)$ satisfying $\omega(x)=\omega(y)$ thus $\omega(x)$ is transitive.

It remains to prove that $X$ has only finitely many sectional-hyperbolic transitive Lyapunov stable sets.
Suppose by absurd that there is an infinite sequence $A_k$ of
sectional-hyperbolic transitive Lyapunov stable sets.
Clearly the members in this sequence must be disjoint, so, since there are finitely many singularities,
we can assume that none of them have singularities.
It follows that all these sets are hyperbolic and then they are all nontrivial hyperbolic attractors of $X$.
In particular, every $A_k$ has dense periodic orbits by the Anosov closing lemma.
We can assume that there is $1\leq i\leq n-2$ such that each $\Lambda_k$ belong to
$Cl(Per_i(X))$.
Define
$$
\Lambda=Cl\left(\bigcup_{k\in \mathbb{N}}A_k\right).
$$
Notice that $\Lambda$ contains infinitely many attractors (the $A_k$'s say).
Moreover, $\Lambda$ is a strongly homogeneous set of index $Ind(\Lambda)=i$
with dense periodic orbits (since each $A_k$ does).

Let us prove that $\Lambda$ satisfies (\ref{eql1}).
Indeed, suppose by contradiction that it does not, i.e.,
there is $\sigma\in Sing_X(\Lambda)$ such that
$I(\sigma)\leq Ind(\Lambda)$.
By Lemma \ref{l4} there is
a dominated splitting $\hat{E}^u_\sigma=E^{uu}_\sigma\oplus E^c_\sigma$
for which the associated strong unstable manifold
$W^{uu}(\sigma)$ satisfies
$\Lambda\cap W^{uu}(\sigma)=\{\sigma\}$.

Take a sequence $x_k\in A_k$ converging to some point $x\in W^s(\sigma)\setminus \{\sigma\}$.
By Corollary 1 p. 949 in \cite{gwz} there is a dominated splitting $D=\Delta^s\oplus \Delta^u$ for the
linear Poincar\'e flow $\psi_t$ which, in virtue of Lemma 2.2 in \cite{gwz},
satisfies $\lim_{t\to\infty}\psi_t(\Delta^u_x)= E^{uu}_\sigma$.
Using exponential maps we can take
a codimension one submanifold $\Sigma$ orthogonal to $X$
of the form $\Sigma=\Delta^s_x(\delta)\times \Delta^u_x(\delta)$
where $\Delta^*_x(\delta)$ indicates the closed $\delta$-ball
around $x$ in $\Delta^*_x$ ($*=s,u$).
Since $\psi_t(\Delta^u_x)\to E^{uu}_\sigma$ as $t\to\infty$
we can assume by replacing $x$ by $X_t(x)$ with $t>0$ large if necessary
that $\Delta_x^u(\delta)$ is almost parallel to $E^{uu}_\sigma$.
In particular, since $\Lambda\cap W^{uu}(\sigma)=\{\sigma\}$,
one has
$(\partial\Delta_x^s(\delta)\times \Delta^u_x(\delta))\cap \Lambda=\emptyset$
where $\partial(\cdot)$ indicates the boundary operation.
Since both $\partial\Delta_x^s(\delta)\times \Delta^u_x(\delta)$ and $\Lambda$
are closed we can arrange a neighborhood $U$ of
$\partial\Delta_x^s(\delta)\times \Delta^u_x(\delta)$ in $\Sigma$
such that $U\cap \Lambda=\emptyset$.

Now we consider $k$ large in a way that $x_k$ is close to $x$.
Replacing $x_k$ by $X_t(x_k)$ with suitable $t$ we can assume that
$x_k\in \Sigma$.
Since $x_k\in A_k$ and $A_k$ is a hyperbolic attractor we
can consider the intersection $S=W^u(x_k)\cap \Sigma$ of the unstable manifold
of $x_k$ and $\Sigma$.
It turns out that $S$ is
the graph of a $C^1$ map
$S:\Delta^u_x(\rho)\to \Delta^s_x(\delta)$ for some
$0<\rho\leq \delta$ whose tangent space $T_yS$ is almost parallel to $\Delta^u_x$.
We assert that $\rho=\delta$.
Otherwise, it would exist some boundary point $z\in \partial S$ in the interior of
$\Sigma$. Since $A_k$ is a hyperbolic set and $z\in A_k$ we could consider
as in \cite{m'} the unstable manifold $W^u(z)$ which will overlap $W^u(x)$.
Since $z\in Int(W^u(z))$ and $W^u(z)\subset A_k$ (for $\Lambda_k$ is an
attractor) we would obtain that $z$ is not a boundary point of $S$, a contradiction which
proves the assertion. It follows from the assertion that
$A_k$ (and so $\Lambda$) would intersect $U$ which is absurd since
$U\cap \Lambda=\emptyset$.
Thus (\ref{eql1}) holds.

Then, Lemma \ref{l3} implies that $\Lambda$ is sectional-hyperbolic for $X$
and so $\Lambda$ has finitely many attractors by Lemma \ref{l5}. But, as we already observed, $\Lambda$ contains
infinitely many attractors so we obtain a contradiction.
This contradiction proves
the finiteness of sectional-hyperbolic transitive Lyapunov stable
sets for $X$ thus ending the proof of the theorem.
\end{proof}

\end{document}